\documentclass[leqno,12pt]{amsart}
\makeatletter
\@namedef{subjclassname@2020}{\textup{2020} Mathematics Subject Classification}
\makeatother
\usepackage{amssymb,enumerate,url}
\theoremstyle{definition}
\newtheorem{example}{\indent Example}

\begin{document}

\title[An example of a domain in a complex flag manifold]{An example of a domain in a complex flag manifold}

\author[N.~Boumuki]{Nobutaka Boumuki}

\address{Division of Mathematical Sciences, Faculty of Science and Technology\endgraf
Oita University, 700 Dannoharu, Oita-shi, Oita 870-1192, JAPAN}
\email{boumuki@oita-u.ac.jp}

\date{}

\subjclass[2020]{Primary 32M10; Secondary 17B20}

\keywords{domain, complex flag manifold, degenerate, constant holomorphic function.}

\maketitle

\begin{abstract}
   This paper shows an example of a connected open neighborhood of a compact connected complex submanifold in a complex flag manifold. 
\end{abstract}

\section{A counter-example}\label{sec-1}
   We shall give an example of a complex flag manifold $X=G_\mathbb{C}/B$ of simply connected complex semisimple Lie group $G_\mathbb{C}$, a closed connected complex (Lie) subgroup $Q\subsetneq G_\mathbb{C}$, a compact connected complex submanifold $Z$ in $X$, and a domain $D$ in $X$ such that
\begin{enumerate}[(c1)]
\item
   $B\subset Q$,
\item
   $Z$ is contained in a fiber of the natural projection $\Pr:G_\mathbb{C}/B\to G_\mathbb{C}/Q$, $gB\mapsto gQ$,   
\item
   $Z\subset D$, and
\item
   $\mathcal{O}(D)\cong\mathbb{C}$, i.e., all holomorphic functions on $D$ are constant.   
\end{enumerate}
   This (c2) and Definition 1 in \cite[p.92]{K} imply that $Z$ is {\it degenerate} in $X$, and hence (c2) and (c4) tell us that the two conditions (1) and (16) in Theorem 2 \cite[pp.92--93]{K} are not equivalent to each other.
   That has effects on some papers (e.g.\ \cite{B2}, \cite{HHL}, \cite{Hu}, \cite{T}, etc).
\par

   Using the notation in Helgason \cite[pp.444--446]{H}, we give

\begin{example}\label{ex-1}
   Let $G_\mathbb{C}:=Sp(2,\mathbb{C})$, $G:=Sp(1,1)$, and $\frak{g}_\mathbb{C}:=\frak{sp}(2,\mathbb{C})$.
   Define two subalgebras $\frak{b},\frak{q}\subset\frak{g}_\mathbb{C}$ and a compact connected subgroup $S\subset G$ by
\[
\begin{split}
&  \frak{b}
   :=\!\left\{\!\begin{array}{c|l}
   \begin{pmatrix}
   z_1 & 0   & 0    & 0    \\
   w_1 & z_2 & 0    & 0    \\
   w_2 & w_3 & -z_1 & -w_1 \\
   w_3 & w_4 & 0    & -z_2 \\
   \end{pmatrix}
   & z_1,z_2,w_1,w_2,w_3,w_4\in\mathbb{C}
              \end{array}\!\right\}\!,\\
&  \frak{q}
   :=\!\left\{\!\begin{array}{c|l}
   \begin{pmatrix}
   z_1 & 0   & 0    & 0    \\
   w_1 & z_2 & 0    & z_3  \\
   w_2 & w_3 & -z_1 & -w_1 \\
   w_3 & z_4 & 0    & -z_2 \\
   \end{pmatrix}
   & z_1,z_2,z_3,z_4,w_1,w_2,w_3\in\mathbb{C}
              \end{array}\!\right\}\!;\\
&  S
   :=\!\left\{\!\begin{array}{c|l}
   \begin{pmatrix}
   x & 0   & 0   & 0    \\
   0 & y_1 & 0   & y_2  \\
   0 & 0   & 1/x & 0 \\
   0 & y_3 & 0   & y_4 \\
   \end{pmatrix}
   & x\in U(1), \begin{pmatrix} y_1 & y_2 \\ y_3 & y_4 \end{pmatrix}\!\in Sp(1)
              \end{array}\!\right\} 
   \,(\cong U(1)\times Sp(1)).           
\end{split}
\]
   Denote by $B$ and $Q$ the normalizers of $\frak{b}$ and $\frak{q}$ in $G_\mathbb{C}$, respectively, and note that $B$ is a Borel subgroup of $G_\mathbb{C}$ (because $\frak{b}$ coincides with the direct sum $\frak{h}\oplus\bigoplus_{\alpha<0}\frak{g}^\alpha$ of the Cartan subalgebra $\frak{h}$ and all the negative-root subspaces $\frak{g}^\alpha\subset\frak{g}_\mathbb{C}$ corresponding to \cite[pp.189--190, The algebra $\frak{c}_n=\frak{sp}(n,\mathbb{C})$ in case $n=2$]{H}) and that $Q$ is a parabolic subgroup of $G_\mathbb{C}$ satisfying (c1) $B\subset Q$.
   In addition, denote by $Z$ (resp.\ $D$) the orbit of the subgroup $S\subset G_\mathbb{C}$ (resp.\ $G\subset G_\mathbb{C}$) through the origin $o\in G_\mathbb{C}/B$, where we remark that $G_\mathbb{C}$ holomorphically acts on $G_\mathbb{C}/B$ in a natural way.
   Then, (c3) $Z\subset D$ holds, and one obtains
\begin{enumerate}[(i)]
\item
   $S\cap B
   =G\cap B
   =\!\left\{\!\begin{array}{c|l}
   \begin{pmatrix}
   x & 0 & 0   & 0   \\
   0 & y & 0   & 0   \\
   0 & 0 & 1/x & 0   \\
   0 & 0 & 0   & 1/y \\
   \end{pmatrix}
   & x,y\in U(1)
              \end{array}\!\right\}\!
   \cong U(1)\times U(1)$,
\item
   $Z=SB/B\cong S/(S\cap B)\cong(U(1)\times Sp(1))/(U(1)\times U(1))$, $\dim_\mathbb{R}Z=2$,
\item
   $D=GB/B\cong G/(G\cap B)\cong Sp(1,1)/(U(1)\times U(1))$, $\dim_\mathbb{R}D=8=\dim_\mathbb{R}G_\mathbb{C}/B$,   
\item
   $SB\subset Q$,
\item
   $Z=SB/B\subset Q/B=\Pr^{-1}(o)$
\end{enumerate}
by a direct computation.
   This (v) implies that (c2) holds, and the (iii) yields (c4) $\mathcal{O}(D)\cong\mathbb{C}$ because $D$ is a flag domain in $G_\mathbb{C}/B$ and the Lie group $G=Sp(1,1)$ is not of Hermitian type (e.g.\ the proof of Lemma 3.4 in \cite[p.5]{B}).
\end{example}

\end{document}